\documentclass[11pt,doublespace]{amsart}
\usepackage{amscd,amssymb,verbatim}
\usepackage{amssymb,amsmath,amsthm,epsfig}

\def\cnum#1{\bigcirc\kern -8pt#1}
\newcommand{\n}{\noindent}

\theoremstyle{plain}
\newtheorem{thm}{Theorem}[section]

\newtheorem{prop}[thm]{Proposition}
\newtheorem{lem}[thm]{Lemma}
\newtheorem{defn}[thm]{Definition}
\newtheorem{rem}[thm]{Remark}
\newtheorem{question}[thm]{Problem}

\newcommand{\e}{\varepsilon}
\newcommand{\T}{\mathcal T}
\newcommand{\B}{\mathcal B}

\newcommand{\cL}{\mathcal L}

\newcommand{\Subs}{\text{Subs}}

\newcommand{\N}{\mathbb{N}}

\newcommand{\Q}{\mathbb{Q}}
\newcommand{\cS}{\mathcal S}
\newcommand{\ssh}{\mathcal{SS}}

\newcommand{\mnorm}[1]{\left\vert\kern-0.9pt\left\vert\kern-0.9pt\left\vert #1
    \right\vert\kern-0.9pt\right\vert\kern-0.9pt\right\vert}

\begin{document}

\title{An Ordinal Index on the Space of Strictly Singular Operators}

\author{Kevin Beanland}
\date{}

\begin{abstract}
Using the notion of $\cS_\xi$-strictly singular operator introduced by Androulakis, Dodos, Sirotkin and Troitsky, we define an ordinal index on the subspace of strictly singular operators between two separable Banach spaces.  In our main result, we provide a sufficient condition implying that this index is bounded by $\omega_1$. In particular, we apply this result to study operators on totally incomparable spaces, hereditarily indecomposable spaces and spaces with few operators.
\end{abstract}

\maketitle

\markboth{Kevin Beanland}{Ordinal Index on the Space of Strictly Singular Operators}

\section{Introduction} \label{sec1}

Dating all the way back to Banach's book \cite{Ba} ordinal indices have been used in Banach space theory to study and measure the complexity of structures present in separable Banach spaces.  Presently, the two most famous examples are the Szlenk and Bourgain indexes \cite{Bo,Sz} (see \cite{Od} for a excellent exposition of both).  More recently, the Schreier families $(\cS_\xi)_{\xi <\omega_1}$ (introduced in \cite{AA}) have been used to index subclasses of classes of separable Banach spaces and classes of operators between separable Banach spaces \cite{AB2,ADST}.  

In a typical example, we consider a class of separable Banach spaces (or subspace of operators) $\mathcal{A}$ all sharing some property $(P)$. For each $1 \leq \xi < \omega_1$, we then define a subclass $\mathcal{A}_\xi$ as those spaces in $\mathcal{A}$ having some stronger property $(P_\xi)$ (whose definition may involve the Schreier family $\cS_\xi$). The goal is to retrieve more information about property $(P)$ by showing that for every $X \in \mathcal{A}$ there is an $1 \leq \xi < \omega_1$ such that $X\in \mathcal{A}_\xi$.

For Banach spaces $X$ and $Y$ a bounded linear operator $T : X \to Y$ is called {\em strictly singular} if the restriction of $T$ to any infinite dimensional subspace is not an isomorphism.  Let $\cL(X,Y)$, $\mathcal{K}(X,Y)$ and $\ssh(X,Y)$ denote the spaces of bounded linear operators, compact operators and strictly singular operators respectively.  In \cite{ADST}, for each $1 \leq \xi < \omega_1$, they define a subset
of $\ssh(X,Y)$ which they call the {\em $\cS_\xi$-strictly singular} operators, denote this set $\ssh_\xi(X,Y)$.  They prove that $\ssh_\xi(X,Y) \subset \ssh_\zeta(X,Y)$ for $\xi< \zeta$ and that for $X$ and $Y$ are separable, every strictly singular operator is {\em $\cS_\xi$-strictly singular} for some $1 \leq \xi < \omega_1$. In other words, for $X$ and $Y$ separable,
$$\bigcup_{\xi < \omega_1} \ssh_\xi(X,Y) =\ssh(X,Y).$$

In this paper we use this definition to introduce the following ordinal index on $\ssh(X,Y)$: Let $X$ and $Y$ be separable Banach spaces. Define,
  
$$i(\ssh(X,Y))= \sup_{T \in \ssh(X,Y)} \inf \{\xi :  T \in \ssh_\xi(X,Y) \}.$$ 

\n Notice that $i(\ssh(X,Y))<\omega_1$ if and only if $\ssh_\xi(X,Y) = \ssh(X,Y)$ for some $1\leq \xi < \omega_1$.

As our main result (Theorem \ref{main}), we provide a simple descriptive set theoretic condition on the space $\ssh(X,Y)$ implying that this index is bounded.  One consequence of this theorem is that if $X$ and $Y$ are totally incomparable (i.e. no subspace of $X$ is isomorphic to a subspace of $Y$) then there is a $\xi<\omega_1$ such that $\ssh_\xi(X,Y)=\ssh(X,Y)$.   We use this theorem to obtain results relating to operators on hereditarily indecomposable Banach spaces, spaces with `few operators.' We also show that for Pelczynski's universal space $U$ \cite{Pe}, $i(\ssh(U))=\omega_1$.

Let us note that this index differs from the Bourgain index in the sense that it is measuring the `richness' of the space $\ssh(X,Y)$ and not necessarily the embeddability of $X$ into $Y$ or the isomorphic comparability of their subspaces.  To illustrate this, recall that $\ell_2$ is isomorphic to all of its infinite dimensional subspace and yet $i(\ssh(\ell_2))=1$ since $\mathcal{K}(\ell_2)= \ssh(\ell_2)$.

In an effort to make this paper accessible to those not familiar with techniques in descriptive set theory we have included, in Section 2, a review of basic definitions in this area that we need to prove our main result. For readers interested in seeing many other applications of descriptive set theoretic methods to Banach space theory, we warmly recommend the paper of Argyros and Dodos \cite{ADo}, the handbook article \cite{AGR} and the recent remarkable work of Dodos \cite{D}.  For the most part, our notation and terminology is standard and can be found in \cite{Ke,LT}.\\ 

\n {\bf Acknowledgments.} We warmly thank both George Androulakis and Spiros Argyros for their guidance during the preparation in this article and interest in this work.  In particular, we thank Prof. Argyros for his hospitality during the authors' visit to NTUA during June 2008 when much of this work was conducted.

\section{Schreier families, Polish spaces and well-founded trees} \label{sec2}

\subsection{Schreier Families}

We begin by recalling the definition of the Schreier families introduced by Alspach and Argyros in \cite{AA}.

If $A$ and $B$ are two finite subsets of $\N$, write $A \leq B$ if $\max A \leq \min B$.  For $n \in N$ and $A \subset \N$, $n \leq A$ if $\{n\} \leq \min A$.  By convention, let $\emptyset < A$ and $A < \emptyset$ for all finite subsets $A$ of $\N$.  For any ordinal number $0 \leq \xi < \omega_1$, the Schreier familiy $\cS_\xi$ is collection of finite subsets of $\N$ defined by the following transfinite recursive process:  

Let, $\cS_0=\{ \{ n \} : n \in \mathbb{N} \} \} \cup \{\emptyset\}$.  Assuming $\xi$ is a successor ordinal and $\cS_\zeta$ has been defined for $\zeta +1= \xi$ let,

$$\cS_{\xi} = \left \{ \bigcup_{i=1}^n F_i : n \geq 1, n \leq F_1 < \cdots < F_n, \text{ and } F_i \in \cS_\zeta \text{ for } 1 \leq i \leq n \right \} \cup \{ \emptyset \}.$$

\noindent For each $n \in \N$ and $\xi < \omega_1$, let $\cS_\xi([n,\infty)) = \{ F \in \cS_\xi : n \leq F\}$.  If $\xi < \omega_1$ is a limit ordinal and $\cS_\zeta$ has been defined for all $\zeta < \xi$ fix an increasing sequence $(\xi_n)_{n=1}^\infty$ such that $\lim_{n \to \infty} \xi_n= \xi$ and define,

$$\cS_{\xi} =\bigcup_{n=1}^\infty \cS_{\xi_n}([n,\infty)).$$

\subsection{Trees on $\N$} Each Schreier family is an example of a {\em tree on $\N$}.  A tree on $\N$ is a collection of finite sequences of $\N$ closed under the partial order of initial segment inclusion.  Let $\N^{<\N}$ denote the collection of all finite sequences of $\N$.  Let $\mathcal{T}r$ denote the set of all trees on $\N$ and $2^{{\N}^{<\N}}$ be the set of all functions from $\N^{<\N}$ to the two point space. $\mathcal{T}r$ is identified with a subset of $2^{{\N}^{<\N}}$ through $\mathcal{T} \mapsto \chi_\mathcal{T}$ (where $ \chi_\mathcal{T}$ is the characteristic function on $\mathcal{T}$).  It is easily seen that $\{\chi_\mathcal{T} : \mathcal{T} \in \mathcal{T}r\}$ is closed in $2^{{\N}^{<\N}}$ when $2^{{\N}^{<\N}}$ is endowed with the product of discrete topologies.  $\mathcal{T} \in \mathcal{T}r$ is called well-founded if  there does not exist a infinite sequence of natural numbers $(\ell_i)_{i=1}^\infty$ such that $(\ell_i)_{i=1}^n \in \mathcal{T}r$ for every $n \in \N$.  Let $\mathcal{WF}$ denote the subset of $\mathcal{T}r$ consisting of all well-founded trees.

We now define the order of a tree.  Let $\prec$ be the partial order of (strict) end-extension.  For every $\mathcal{T} \in \mathcal{T}r$ let 
$$\mathcal{T}' = \{ s \in \mathcal{T} : \text{ there exists } t \in \mathcal{T} \text{ with } s \prec t \}.$$
Observe that $\mathcal{T}' \in \mathcal{T}r$.  By transfinite recursion, for every $\mathcal{T} \in \mathcal{T}r$ we define $(\mathcal{T}^{(\xi)})_{\xi< \omega_1}$ as follows:   

$$\mathcal{T}^{(0)} = \mathcal{T},~ \mathcal{T}^{(\xi + 1)} = (\mathcal{T}^{(\xi)})' \text{ and } \mathcal{T}^{(\lambda)} = \bigcap_{\xi < \lambda} \mathcal{T}^{(\xi)},$$

\n whenever $\lambda$ is a limit ordinal.  $\mathcal{T} \in \mathcal{WF}$ if and only if the sequence $(\mathcal{T}^{(\xi)})_{\xi<\omega_1}$ is eventually empty.  For every $\mathcal{T} \in \mathcal{WF}$ the order  of $\mathcal{T}$, denoted $o(\mathcal{T})$, is defined to be the least countable ordinal $\xi$ such that $\mathcal{T}^{(\xi)} = \emptyset$.   For each $1 \leq \xi < \omega_1$, $\cS_\xi \subset \mathcal{WF}$.  The next proposition computes the height of the trees $o(\cS_\xi)$.

\begin{prop}
For each $1 \leq \xi < \omega_1$, $o(\cS_\xi) = \omega^\xi$. \label{Schreierheight}
\end{prop}

\subsection{Polish spaces, Borel and analytic subsets} The main objects of study in descriptive set theory are separable metrizable spaces or {\em Polish} spaces.  The canonical example of a Polish space is the so-called Baire space $\N^\N$ equipped with the product of discrete topologies.  The Borel subsets of a Polish spaces are those sets obtained by taking countable intersections, unions and complements of the open sets.   Given a space $P$ and its topology $\mathfrak{T}$, denote by $B(\mathfrak{T})$ the  Borel subsets generated by $\mathfrak{T}$.  A subset $\mathcal{A}$ of a Polish space $P$ is {\em analytic} if it is the continuous image of the Baire space and {\em coanalytic} if it is the complement of such a set. We collect two useful facts about these sets are in the following remark.
\begin{rem} Let $(\mathcal{A}_n)_n$ be countable collection of analytic subsets of a Polish space $P$.
\begin{enumerate}
\item $\cap_n \mathcal{A}_n$ and $\cup_n \mathcal{A}_n$ are both analytic.
\item A set is Borel if and only if it is both analytic and coanalytic.  
\end{enumerate}\label{analyticrem}
\end{rem}
A function between Polish spaces is Borel if the inverse image of any Borel subset is Borel. The next proposition characterizes of the definition of analytic.  This characterization will be used repeatedly throughout this paper.

\begin{prop}
For a Polish space $P$, a subset $\mathcal{A}$ of $P$ is analytic if and only if there is a Polish space $S$ and Borel map $\phi:S \to P$ such that $\phi(S)=\mathcal{A}$. \label{analyticdefn}
\end{prop}

Since $\mathcal{T}r$ is a closed subset of the Polish space $2^{{\N}^{<\N}}$ it is itself a Polish space in the relative topology.  The next theorem is known as the {\em Boundedness Theorem for Well-Founded Trees} \cite[Theorem 31.2]{Ke}.

\begin{thm} \label{boundedness}
Let $\mathcal{A}$ be an analytic collection of well founded trees on $\N$. Then 
$$\sup\{o(\T): \T \in \mathcal{A} \} < \omega_1.$$
\end{thm}
The following classical fact allows one to regard any Borel subset of a Polish space as a Polish space in its own right \cite[Corollary 13.4]{Ke}.

\begin{thm}\label{borel}
Let $P$ be a Polish space with the topology $\mathfrak{T}$ and $S \in B(\mathfrak{T})$.  There exists a finer Polish topology $\mathfrak{U}$ on $P$ such that,
\begin{itemize}
\item[(1)] $S$ is clopen in $\mathfrak{U}$;
\item[(2)] $B(\mathfrak{T})=B(\mathfrak{U})$.
\end{itemize}
Moreover, $S$ with the relative topology of $\mathfrak{U}$ is a Polish space. \label{borelpolish}
\end{thm}

Let $X$ and $Y$ be separable Banach spaces. The following topological spaces are Polish spaces in the topologies of coordinatewise and pointwise convergence respectively,
$$X^\N=\{(x_n)_{n=1}^\infty: x_n \in X \text{ for all } n \in \N \} \quad and \quad \B(\cL(X,Y))=\{T \in \cL(X,Y): \|T\| \leq 1 \}.$$
 \n See \cite{Ke} for details regarding $\B(\cL(X,Y))$.  The topology of pointwise convergence on 
 $\B(\cL(X,Y))$  is called the {\em strong-operator} topology.

For any Banach space $X$, let $\mathfrak{B}_X$ denote the set of all normalized basic sequences in $X$.  For $k \in \N$, let $\mathfrak{B}_X^k$ be the subset of $\mathfrak{B}_X$ with basis constant less that $k$.  Since $\mathfrak{B}_X = \cup_{k=1}^\infty \mathfrak{B}_X^k$ and each $\mathfrak{B}_X^k$ is closed in $X^\N$, $\mathfrak{B}_X$ is a Borel subset of $X^\N$ (in particular it is  $F_\sigma$).  Invoking Proposition \ref{borelpolish}, we know that $\mathfrak{B}_X$ is a Polish space whose Borel sets coincide with the relative Borel subsets of a $X^\N$.  

Let $X$ be a separable Banach space. Let $F(X)$ denote the set of all closed subsets of $X$ and $\Subs(X)$ denote the subset of all closed infinite dimensional subspaces of $X$.  Let $F(X)$ be endowed with the  $\sigma$-algebra $\Sigma$ generated by

$$\{ F \in F(X) : F \cap U \not= \emptyset \}$$

\n where $U$ ranges over all non-empty open subsets of $X$. $(F(X),\Sigma)$ is called the {\em Effros-Borel} space of $X$.  It is well-known that there is a Polish topology $\mathfrak{T}$ on $F(X)$ such that $B(\mathfrak{T})=\Sigma$.  Since $\Subs(X) \in \B(\mathfrak{T})$ it may be regarded as a Polish space it is own right.

Finally, we will be using the following consequence of the famous Kuratowski-Ryll-Nardzewski selection theorem (see \cite[Theorem 12.13]{Ke}).

\begin{thm}
Let $X$ be a separable Banach space. There is a sequence of maps $(d_\ell)_{\ell=1}^\infty$ such that,
\begin{enumerate}
\item$d_\ell:Subs(X) \to X$ and $d_\ell$ is Borel for each $\ell \in \N$;  
\item For each $Y \in \Subs(X)$, $(d_\ell(Y))_{\ell=1}^\infty$ is norm dense in $S_Y$.
\end{enumerate}\label{KRN}
\end{thm}

\section{Main Results}\label{sec3}

In our first proposition, we give an upper bound on the complexity $\mathcal{B}(\ssh(X,Y)):= \ssh(X,Y)\cap\mathcal{B}(\cL(X,Y))$ as a subset of $\B(\cL(X,Y))$ in the strong operator topology. As we will observe later, for certain classes of Banach spaces this bound is sharp and for others it can be reduced.

\begin{prop}
Let $X$ and $Y$ be separable Banach spaces.  Then $\mathcal{B}(\ssh(X,Y))$ is a coanalytic subset of $\mathcal{B}(\cL(X,Y))$ in the strong operator topology. \label{coanalytic}
\end{prop}

\begin{proof}
By Proposition \ref{KRN} there is a sequence of Borel maps $(d_\ell)_{\ell=1}^\infty$ from $\Subs(X)$ to $X$ such 
that for each $Z \in \Subs(X)$, $(d_\ell(Z))_{\ell=1}^\infty$ is norm dense in $S_Z$. Observe that,
\begin{equation*}
\begin{split}
R \in \B(\cL(X,Y))\setminus\B(\ssh(X,Y))  \iff& \exists~ Z \in \Subs(X),~ \exists~ n \in \N ~ \text{ such that }\\
& \forall \ell \in \N,~ \|R (d_\ell(Z)) \| \geq 1/n. \\
\end{split}
\end{equation*}

The fact that there is a existential quantifier over a Borel set ($\Subs(X)$) and the remaining quantifiers are over countable sets indicates that $\B(\cL(X,Y))\setminus\B(\ssh(X,Y))$ is a analytic set.   However, for the convenience of readers not familiar with descriptive set theory, we give a more detailed argument.  Let,

$$\mathcal{B}_n= \{(R,Z)\in \B(\cL(X,Y))\times \Subs(X) : \|R( d_\ell(Z)) \| \geq 1/n  \text{ for all } \ell \in \N \}.$$

\n Let $\pi_1$ be the projection of $\B(\cL(X,Y))\times \Subs(X)$ onto $\B(\cL(X,Y))$.  Then 
$$\B(\cL(X,Y))\setminus \B(\ssh(X,Y))= \bigcup_{n=1}^\infty \pi_1(\mathcal{B}_n).$$  

\n Invoking Remark \ref{analyticrem} and Proposition \ref{analyticdefn}, it suffices to show that $\mathcal{B}_n$ is Borel for each $n \in \N$.  For each $\ell \in \N$ define $H_\ell:\B(\cL(X,Y))\times \Subs(X) \to \mathbb{R}$ by $H_\ell(R,Z)=\|R( d_\ell(Z)) \|$.  By the continuity of the norm and the Borelness of the map $d_\ell$, the map $H_\ell$ is Borel.  Since $\mathcal{B}_n=\cap_{\ell=1}^\infty H_\ell^{-1}[1/n,\infty)$ is Borel, the claim follows.
\end{proof}

The following remark gives an simple characterization of strict singularity in terms of the behavior of the operator on finitely supported vectors of normalized basic sequences and, in turn, motivates the definition $\cS_\xi$-strictly singular.  For $(x_n)_n \in \mathfrak{B}_X$ and $A \subset N$, let $[x_n]_{n\in A}$ denote the closed linear span of the vectors $(x_n)_{n\in A}$.  

\begin{rem}
$T \in \ssh(X,Y)$ if and only if for every $\e >0$ and $(x_n)_n \in \mathfrak{B}_X$ there is a finite subset $F$ of $\N$ 
and a vector $x \in [x_n]_{n \in F} \setminus \{0\}$ such that $\| T x\| < \e \|x\|$.  
\end{rem}

From this characterization it is natural to define a subset of the strictly singular operators in the following way:  Let $\mathcal{A}$ be any collection of finite subsets $\N$.  By replacing the condition `there is a finite subset $F$ of $\N$' in the above remark with `$F \in \mathcal{A}$' we arrive at a new class of operators which is appropriate referred to as $\mathcal{A}$-strictly singular.  By definition, every $\mathcal{A}$-strictly singular is strictly singular.  It is in this way that the $\cS_\xi$-strictly singular operators are defined.  

\begin{defn} Let $1 \leq \xi < \omega_1$.  
$T \in \cL(X,Y)$ is $\cS_\xi$-strictly singular, written $T \in \ssh_\xi(X,Y)$, if and only if for every $\e >0$ and $(x_n)_n \in \mathfrak{B}_X$ there is a $F \in \cS_\xi$
and a vector $x \in [x_n]_{n \in F} \setminus \{0\}$ such that $\| T x\| < \e \|x\|$.  
\end{defn}

We refer the reader to \cite{ADST} for a detailed account of the many interesting properities of these operators. The next result is the main result of this paper.  The proof follows the outline of the proof of Theorem 6.5(ii) in \cite{ADST}.  The main difference is the incorporation of the Polish space $\B(\cL(X,Y))$.

\begin{thm}
Let $X$ and $Y$ be separable Banach spaces. Suppose that $\B(\ssh(X,Y))$ is Borel subset of $\B(\mathcal{L}(X,Y))$ in the strong operator topology. Then $i(\ssh(X,Y))<\omega_1$.\label{main}
\end{thm}

Observe that by Proposition \ref{coanalytic} and Remark \ref{analyticrem} (2), we may replace the condition {\em Borel} in the above theorem by {\em analytic}.  Recall that  $i(\ssh(X,Y))<\omega_1$ if and only if there exists a $\xi < \omega_1$ such that $\ssh_\xi(X,Y) = \ssh(X,Y)$ 

\begin{proof}
It suffices to show that $\B(\ssh_\xi(X,Y)) = \B(\ssh(X,Y))$. 

For each, $R \in \B(\ssh(X,Y)), m\in \N$ and $(x_n) \in \mathfrak{B}_X$, define a tree on $\N$ in the following way:  
  $$\mathcal{T}(R,m,(x_n)_n)= \{ ( l_1, \dots , l_n) : \forall (a_i)_i \in \Q^{< \N} , \| R(\sum_{i=1}^n a_i x_{l_i} )\| \geq \frac{1}{m} \| \sum_{i=1}^n a_i x_{l_i} \| \},$$
where $\Q^{< \N}$ denotes the set of all finite sequences of rationals.  It is easy to see that $\mathcal{T}(R,m,(x_n)_n)$ is a tree on $\N$.  

Let us see that $\mathcal{T}(R,m,(x_n)_n) \in \mathcal{WF}$. Supposing not, we find an infinite sequence $(l_i)_{i=1}^\infty$ such that $(l_i)_{i=1}^k \in \mathcal{T}(R,m,(x_n)_n)$ for all $k \in \N$.  It follows that for each $(a_i) \in \Q^{< \N}$,  
$$\| R \sum_{i=1}^\infty a_i x_{l_i} \|\geq \frac{1}{m}  \| \sum _{i=1}^\infty a_i x_{l_i} \| .$$  
This implies that $R|_{ [ x_{l_i}  ]_{i=1}^\infty} $ is an isomorphism, contradicting the fact that $R \in \B(\ssh(X,Y))$.  

We wish to show that the following collection of well founded trees is analytic as a subset of $\mathcal{T}r$: 
$$\mathcal{A} = \{ \mathcal{T} ( R, m, (x_n)_n ) : R \in \B(\ssh(X,Y)), m \in \N, (x_n)_n \in \mathfrak{B}_X\}.$$ 
Assume for the moment that $\mathcal{A}$ is analytic. From Theorem \ref{boundedness} we can find an $\xi < \omega_1$ such that $\sup \{ o(\mathcal{T}) : \mathcal{T} \in \mathcal{A} \} < \xi$.  We claim that $\B(\ssh_\xi(X,Y))= \B(\ssh(X,Y))$.  Suppose, for the sake of contradiction, that there exists $R \in \B(\ssh(X,Y)) \setminus \B(\ssh_\xi(X,Y))$.  By the definition of $\B(\ssh_\xi(X,Y)) $, there exists $(x_n) \in \mathfrak{B}_X$ and $m \in \N$ such that for all $F \in \cS_\xi$ and for all $(a_i) \in \Q^{< \N}$, 
$$\|R \sum_{i \in F} a_i x_i \| \geq \frac{1}{m} \| \sum_{i \in F} a_i x_i \|.$$  
This means that $\cS_\xi \subset \mathcal{T}( R, m, (x_n))$.  Applying Proposition \ref{Schreierheight} yields the following contradiction, 
$$ \xi > o(\mathcal{T} (R, m, (x_n))) \geq o(\cS_\alpha) = \omega^\xi.$$ 

Therefore, if suffices to prove that $\mathcal{A}$ is analytic.  $\mathcal{A}$ is the countable union of the following sets:  For each $m \in \N$ let,
$$\mathcal{A}_m = \{ T ( R, m, (x_n)_n ) : R \in \B(\ssh(X,Y)),(x_n)_n \in \mathfrak{B}_X\}.$$
Invoking Remark \ref{analyticrem}(1), it suffices to show that $\mathcal{A}_m$ is analytic for each $m \in \N$.  By Proposition \ref{analyticdefn} we are tasked with finding a Polish space $P$ and a Borel map $\varphi_m:P \to \mathcal{T}r$ such that $\varphi_m(P)=\mathcal{A}$.  Our Polish space will be $\B(\ssh(X,Y)) \times \mathfrak{B}_X$ and our map is the following:
$$\varphi_m(R, (x_n)_n) = \mathcal{T}(R,m,(x_n)_n).$$
Since $\mathfrak{B}_X$ is a Borel subset of $X^\N$ and we have assumed $\B(\ssh(X,Y))$ is a Borel subset of $\B(\mathcal{L}(X,Y))$ the space $\B(\ssh(X,Y)) \times\mathfrak{B}_X $ is a Polish space whose Borel subsets coincide with the relative Borel subsets of $\B(\cL(X,Y))\times X^\N$ (see Proposition \ref{borelpolish}).

The final step is to show that $\varphi_m$ is a Borel map.  It suffices to show that that the inverse image of a basic open neighborhood of $\mathcal{T}r$ is a Borel subset of $\B(\ssh(X,Y)) \times\mathfrak{B}_X $.  The basic open neighborhoods of $\mathcal{T}r$ take the form $U_F = \{ \mathcal{T} \in \mathcal{T}r : F \in \mathcal{T} \}$ where $F$ is a fixed finite subset of $\N$.  Fix $F$ and observe that,  

\begin{equation*}
\begin{split}
\varphi_m^{-1}(U_F)&= \{ (R, (x_n)) : F \in \mathcal{T} (R, m, (x_n)) \} \\
&= \{ (R, (x_n)) : \text{ for all } (a_i) \in \Q^{< \N},~ \| R \sum_{i \in F} a_i x_i \| \geq \frac{1}{m} \| \sum_{i \in F} a_i x_i \| \}\\
& = \bigcap_{ (a_i) \in \Q^{< \N} } \{ (R, (x_n)) : ~ \| R \sum_{i \in F} a_i x_i \| \geq \frac{1}{m} \| \sum_{i \in F} a_i x_i \| \}
\end{split}
\end{equation*}
\n Now fix $(a_i)_i \in \Q^{< \N}$.  Since,
$$\mathcal{C}_{((a_i),F)}=\{ (R, (x_n)) : ~ \| R \sum_{i \in F} a_i x_i \| \geq \frac{1}{m} \| \sum_{i \in F} a_i x_i \| \}$$
\n is clearly closed in $ \B(\cL(X,Y)) \times X^\N$, it is a Borel subset of $\B(\ssh(X,Y)) \times \mathfrak{B}_X$.  This proves that $\mathcal{A}$ is analytic.
\end{proof}

\section{Applications and Further Research}

\subsection{Totally Incomparable Spaces} The first obvious consequence of Theorem \ref{main} is that for separable spaces $X$ and $Y$, if $\cL(X,Y)=\ssh(X,Y)$, then there is an $\xi <\omega_1$ such that $\ssh_\xi(X,Y)=\ssh(X,Y)$.  In particular this is true when the spaces $X$ and $Y$ are totally incomparable.  

\subsection{Hereditarily Indecomposable Spaces} The second example to illustrate our theorem concerns the class of hereditarily indecomposable Banach spaces. A Banach space $X$ is hereditarily indecomposable (HI) if for any two infinite dimensional subspaces $X$ and $Y$ and $\e>0$ there is a $x \in S_X$ and $y \in S_Y$ such that $\|x-y\|<\e$ ($S_X$ denotes the set of all norm-one vectors in $X$).  For more information regarding these important spaces we refer the reader to the works \cite{ATo,AT,GM}.  

The next proposition shows that whenever the space of strictly singular operators between two separable Banach spaces can be characterized in a specific way, the complexity of $\B(\ssh(X,Y))$ reduces from coanalytic to Borel.  

\begin{prop}
Suppose $X$ and $Y$ are separable Banach spaces such that 
$$R \in \ssh(X,Y) \iff \forall~ \e>0, ~\exists~ Z \in \Subs(X)\text{ such that }\|R|_Z\|<\e.$$ 
Then $\B(\ssh(X,Y))$ is a Borel subset of $\B(\cL(X,Y))$.  Moreover, applying Theorem \ref{main},  $i(\ssh(X,Y))< \omega_1$.  \label{ssequiv}
\end{prop}

\begin{proof} 
 Combining Proposition \ref{coanalytic} and Remark \ref{analyticrem} (2) it suffices to show that $\B(\ssh(X,Y))$ is analytic.
By Proposition \ref{KRN} there is a sequence of Borel maps $(d_\ell)_{\ell=1}^\infty$ from $\Subs(X)$ to $X$ such 
that for each $Y \in \Subs(X)$, $(d_\ell(Y))_{\ell=1}^\infty$ is norm dense in $S_Y$.  For each $n \in \N$ let,
$$\mathcal{A}_n=\{(R,Z)\in \B(\cL(X,Y))\times \Subs(X) : \|R( d_\ell(Z)) \| < 1/n \text{ for all } \ell \in \N \}.$$
\n Let $\pi_1$ be the projection of $\B(\cL(X,Y))\times \Subs(X)$ onto $\B(\cL(X,Y))$.  By assumption we can observe that,
$$R\in \B(\ssh(X,Y)) \iff R \in \bigcap_{n=1}^\infty \pi_1(\mathcal{A}_n).$$
\n Therefore showing $\mathcal{A}_n$ is Borel proves our claim.  Since this follows from arguments similar to those found in the proof of Proposition \ref{coanalytic} and we omit the proof. 
\end{proof}

We will show that every strictly singular operator originating from a separable HI space $X$ and having a separable range space $Y$ satisfies the assumptions of Proposition \ref{ssequiv}. Therefore, in this case, $i(\ssh(X,Y)< \omega_1$.  This result should be compared to \cite[Theorem 6.13]{ADST}.

\begin{lem}
Let $X$ be an HI space and $Y$ be any Banach space.  Then $T\in \ssh(X,Y)$ if and only if for all $\e>0$ there is an infinite dimensional subspace $Z$ of $X$
such that $\|T|_Z\|<\e$.\label{HIborel}
\end{lem}
   
\begin{proof}
The forward implication follows from a well-known characterization of strictly singular operators \cite[page 76]{LT}.  

The reverse implication follows an argument found in \cite{ATo}.  Assume without loss of generality that $\|T\|=1$.  Find $\e>0$ and an infinite dimensional subspace $W$ of $X$ such
that $\|Tw \| > \e$ for all $w \in S_W$.  Now find an infinite dimensional subspace $Z$ such that $\| Tz \| < \e/2$ for all $z \in S_Z$.  For $w \in S_W$ and $z \in S_Z$
we have,

$$\|z-w \| \geq \|Tw-Tz\| \geq \e-\frac{\e}{2} =\frac{\e}{2}.$$

\n Therefore, $X$ is not HI; the claim follows.
\end{proof}

\subsection{Spaces admitting few operators}\label{few} A Banach space $X$ is said to admit {\em few operators} if every operator $T$ on $X$ takes the form $T= \lambda I + S$ where $\lambda$ is a scalar and $S$ is a strictly singular operator on $X$.  The first examples of spaces admitting few operators were all HI spaces, however, it has since been shown that such spaces need not be HI.  In particular, there are constructions of spaces exhibiting rich unconditional structure and admitting few operators simultaneously \cite{AM2,AM1}.  The following proposition together with Theorem \ref{ssequiv}, shows that if $X$ admits few operators then $i(\ssh(X))< \omega_1$.

\begin{prop}
Suppose $X$ is separable and admits few operators.  Then $T \in \ssh(X)$ if and only if for all $\e>0$ there is an infinite dimensional subspace $Z$ of $X$ such that $\|T|_Z\|<\e$.
\end{prop}

\begin{proof}
As in the proof of Proposition \ref{HIborel}, the forward direction is automatic.  For the reverse direction, we proceed by contradiction.  Let $T \in \cL(X) \setminus \ssh(X)$ satisfy the assumption. Find a scalar $\lambda$ and a $S \in \ssh(X)$ such that $T=\lambda I + S$.  By assumption, there is an infinite dimensional subspace $Z$ of $X$ such that $\|(\lambda I + S)|_{Z}\|< |\lambda |/4$.    Since $S$ is strictly singular, there is an infinite dimensional subspace $W$ of $Z$ such that $\|S|_W\| < |\lambda|/2$.  Let $w \in S_W$ and observe that,
$$\frac{|\lambda|}{4} > \|(\lambda +S)w\| > |\lambda| - \frac{|\lambda|}{2}  = \frac{|\lambda|}{2} .$$
This contradiction, proves the claim.
\end{proof}

\subsection{A universal space producing an unbounded index} \label{universal} As our final demonstration of Theorem \ref{main}, we show that $i(\ssh(U))=\omega_1$; where $U$ is Pe{\l}czy{\'n}ski's universal space \cite{Pe}.  From Theorem \ref{main}, this will imply that $\B(\ssh(U))$ is a coanalytic non-Borel subset of $\B(\cL(U))$.  

In \cite[Example 2.7]{ADST} they exhibit a collection of spaces $(T_\xi)_{\xi <\omega_1}$ such that $T_{\omega \xi} \subset T_\xi$ (as vector spaces) and for which the identity $i_\xi: T_{\omega \xi} \to T_\xi$ satisfies $i_\xi \in \ssh_{\omega \xi} (T_{\omega \xi} ,T_\xi)\setminus \ssh_{\xi} (T_{\omega \xi} ,T_\xi)$.  Let $P_\xi$ be the projection from $U$ onto $T_\xi$ (assume, for simplicity, that $T_\xi \subset U$).  Then $i_\xi \circ P_\xi \in \ssh(U) \setminus \ssh_\xi(U)$.  

\subsection{Future Research and Concluding Remarks} 

The first natural question that arises is whether the converse
 of Proposition \ref{main} holds.
\begin{question}
Suppose $X$ and $Y$ are separable Banach spaces and there is an $\xi< \omega_1$ such that $\ssh_\xi(X,Y)= \ssh(X,Y)$.  Is $\B(\ssh(X,Y))$ a Borel subset of $\B(\cL(X,Y))$?
\end{question}

This problem was communicated to the author by S. Argyros. We conjecture that it has a positive answer.  

It would be desirable to have a general scheme, perhaps determined by intrinsic properties of the spaces, that could either show how to reduce the complexity of the space $\B(\ssh(X,Y))$ or  exhibit a sequence of operators with properties similar to those in Subsecton \ref{universal}.

\bibliographystyle{plain}
\bibliography{bib_source}

{\footnotesize
\noindent
Department of Mathematics and Applied Mathematics, Virginia Commonwealth University, Richmond, VA 23284, \\
\noindent kbeanland@vcu.edu }

\end{document}